\documentclass[graybox]{svmult}

\usepackage{newtxtext,newtxmath}
\usepackage{siunitx}

\usepackage{helvet}         
\usepackage{courier}        
\usepackage{makeidx}         
\usepackage{graphicx}        
\usepackage[bottom]{footmisc}
\usepackage{url}

\makeindex             

\begin{document}

\title*{Topological Surgery in Nature}
\author{Stathis Antoniou and Sofia Lambropoulou}
\institute{Stathis Antoniou \at Department of Mathematics, National Technical University of Athens, Athens, Greece \\ \email{santoniou@math.ntua.gr}
\and Sofia Lambropoulou  \at Department of Mathematics, National Technical University of Athens, Athens, Greece \\ \email{sofia@math.ntua.gr}}
\maketitle

\abstract{In this paper, we extend the formal definition of topological surgery by introducing new notions in order to model natural phenomena exhibiting it. On the one hand, the common features of the presented natural processes are captured by our schematic models and, on the other hand, our new definitions provide the theoretical setting for examining the topological changes involved in these processes.}

\let\thefootnote\relax\footnotetext{

{\noindent}\textit{2010 Mathematics Subject Classification}: 57R65, 57N12, 57M25, 57M99, 37B99, 92B99

{\noindent}\textit{Keywords}: topological surgery, gluing homeomorphism, three-space, three-sphere, solid topological surgery, embedded  topological surgery, mathematical modeling, natural processes, joining `thread', topological `drilling', continuous process, dynamics, attracting forces, repelling forces, recoupling, reconnection, DNA recombination, chromosomal crossover, magnetic reconnection, cosmic magnetic lines, Falaco solitons, tornadoes, whirls, waterspouts, mitosis, meiosis, gene transfer, necking, black holes}

\setcounter{section}{-1}
\section{Introduction} \label{Intro}
Topological surgery is more than a mathematical definition allowing the creation of new manifolds out of known ones. As we point out in \cite{SSN1},  it is a process triggered by forces that appear in both micro and macro scales. For example,  in dimension 1 topological surgery can be seen in DNA recombination and during the reconnection of cosmic magnetic lines,  while in dimension 2 it happens when genes are transferred in bacteria and during the formation of black holes. Inspired by these natural processes,  in \cite{SSN1}  we enhance the formal definition of topological surgery with the observed {\it forces} and {\it dynamics},  thus turning it into a continuous process.
 
Furthermore,   in \cite{SSN1} we observe that phenomena like tension on soap films or the merging of oil slicks are undergoing $1$-dimensional surgery but they happen on surfaces instead of $1$-manifolds. Similarly,  moving up one dimension,  during the biological process of mitosis and during tornado formation,  $2$-dimensional surgery is taking place on $3$-dimensional manifolds instead of surfaces. Thus,  in order to  fit natural phenomena where the interior of the initial manifold is filled in,  in \cite{SSN1} we extend the formal definition by introducing the notion of {\it solid topological surgery} in both dimensions 1 and 2.

Finally,  in \cite{SSN1} we notice that in some phenomena exhibiting topological surgery,  the ambient space is also involved. For example in dimension 1,  during DNA recombination the initial DNA molecule which is recombined can also be knotted. In other words,  the initial $1$-manifold can be a knot (an embedding of the circle) instead of an abstract circle. Similarly in dimension 2,  the processes of tornado and black hole formation are not confined to the initial manifold and topological surgery is causing (or is caused by) a change in the whole space. We therefore define the notion of {\it embedded topological surgery} which allows us to model these kind of phenomena but also to view all natural phenomena exhibiting topological surgery as happening in $3$-space instead of abstractly.

The notions and ideas presented in this paper can be found in \cite{SSN1, SSN2}. However,  \cite{SSN1, SSN2} contain a much deeper analysis of the dynamical system modeling 2-dimensional surgery,  which is discussed here very briefly.  

The paper is organized as follows: it starts by recalling the formal definitions of surgery in Section~\ref{definitions}. Dynamic topological surgery is then presented in natural processes and introduced as a schematic model in Section~\ref{Dynamics}. Next,  solid surgery is defined in Section~\ref{DSolid} where the dynamical system modeling 2-dimensional surgery is also presented. Finally,  embedded surgery and its differences in dimensions 1 and 2 are discussed in Section~\ref{Embedded}.
 
We hope that the presented definitions and phenomena will broaden our understanding of both the topological changes exhibited in nature and topological surgery itself.

\section{The formal definitions of surgery} \label{definitions}

We recall the following well-known definition:

\begin{definition} \label{surgery} \rm An \textit{m-dimensional n-surgery} is the topological procedure of creating a new $m$-manifold $M'$ out of a given $m$-manifold $M$ by removing a framed $n$-embedding $h:S^n\times D^{m-n}\hookrightarrow  M$ and replacing it with $D^{n+1}\times S^{m-n-1}$,  using the `gluing' homeomorphism $h$ along the common boundary $S^n\times S^{m-n-1}$. Namely,  and denoting surgery by  $\chi$:
\[M' = \chi(M) = \overline{M\setminus h(S^n\times D^{m-n})} \cup_{h|_{S^n\times S^{m-n-1}}} (D^{n+1}\times S^{m-n-1}). \]
 \end{definition}
 
{\noindent}Further,  \rm the \textit{dual m-dimensional $(m-n-1)$-surgery} on $M'$ removes a dual framed $(m-n-1)$-embedding  $g:D^{n+1}\times S^{m-n-1}\hookrightarrow  M'$ such that $g|_{S^n\times S^{m-n-1}}=h^{-1}|_{S^n\times S^{m-n-1}}$,  and replaces it with $S^n\times D^{m-n}$,  using the `gluing' homeomorphism $g$ (or $h^{-1}$) along the common boundary $S^n\times S^{m-n-1}$. That is:
\[M = \chi^{-1}(M') = \overline{M'\setminus g(D^{n+1}\times S^{m-n-1})} \cup_{h^{-1}|_{S^n\times S^{m-n-1}}} (S^n\times D^{m-n}). \]

{\noindent}From the above definition it follows that $M = \chi^{-1}(\chi(M))$ and $n+1 \leq m$. For further reading see \cite{Ra,  PS,  Ro}. We shall now apply the above definition to dimensions 1 and 2.

\subsection{1-dimensional 0-surgery}\label{1D_FormalE}
\begin{figure}[!h]
\begin{center}
\includegraphics[width=11cm]{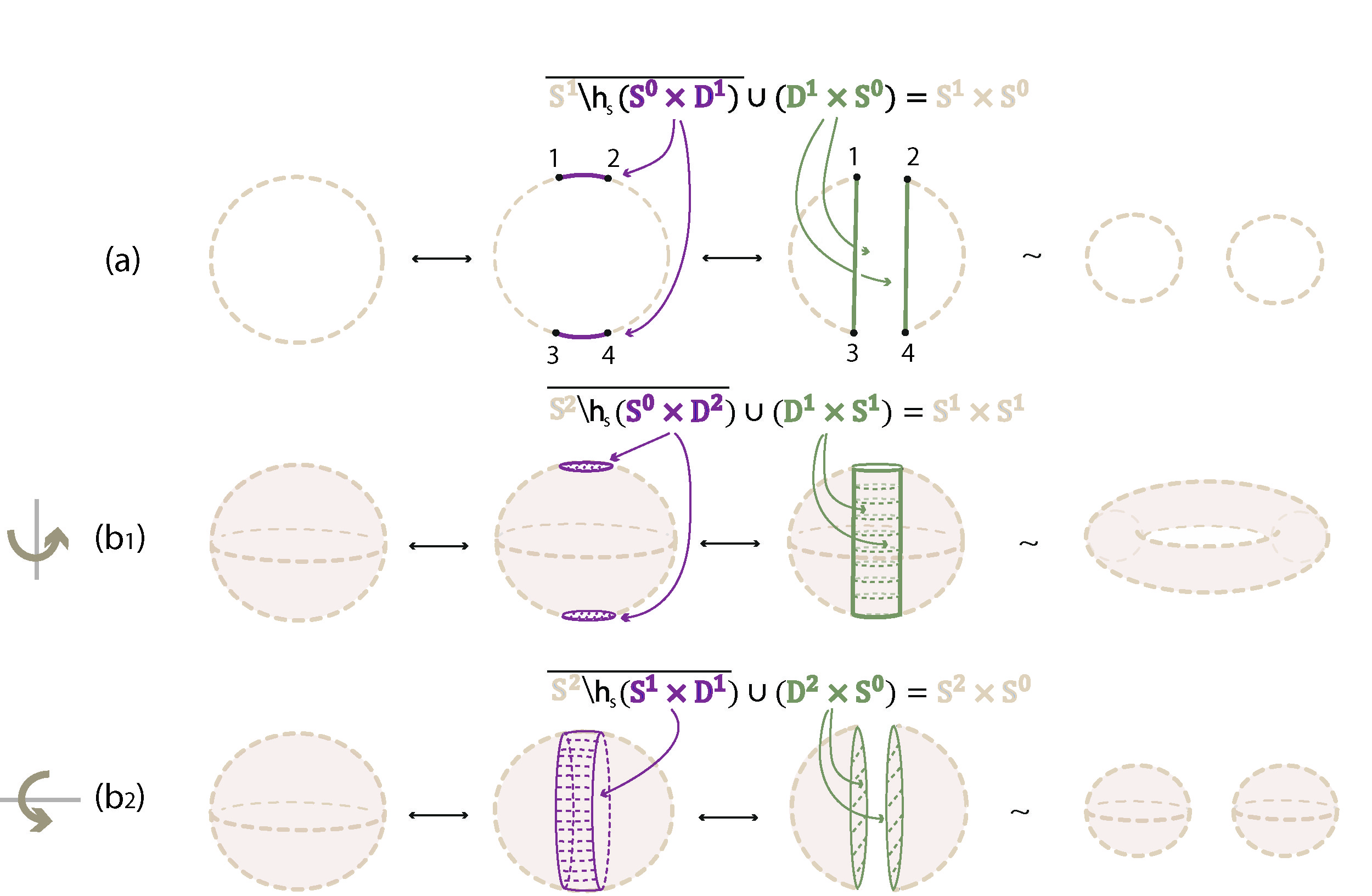}
\caption{Formal \textbf{(a)} 1-dimensional 0-surgery \textbf{(b\textsubscript{1})} 2-dimensional 0-surgery and \textbf{(b\textsubscript{2})} 2-dimensional 1-surgery using the standard embedding $h_s$.}
\label{1_2D_Formal}
\end{center}
\end{figure}

\begin{figure}[!h]
\begin{center}
\includegraphics[width=11cm]{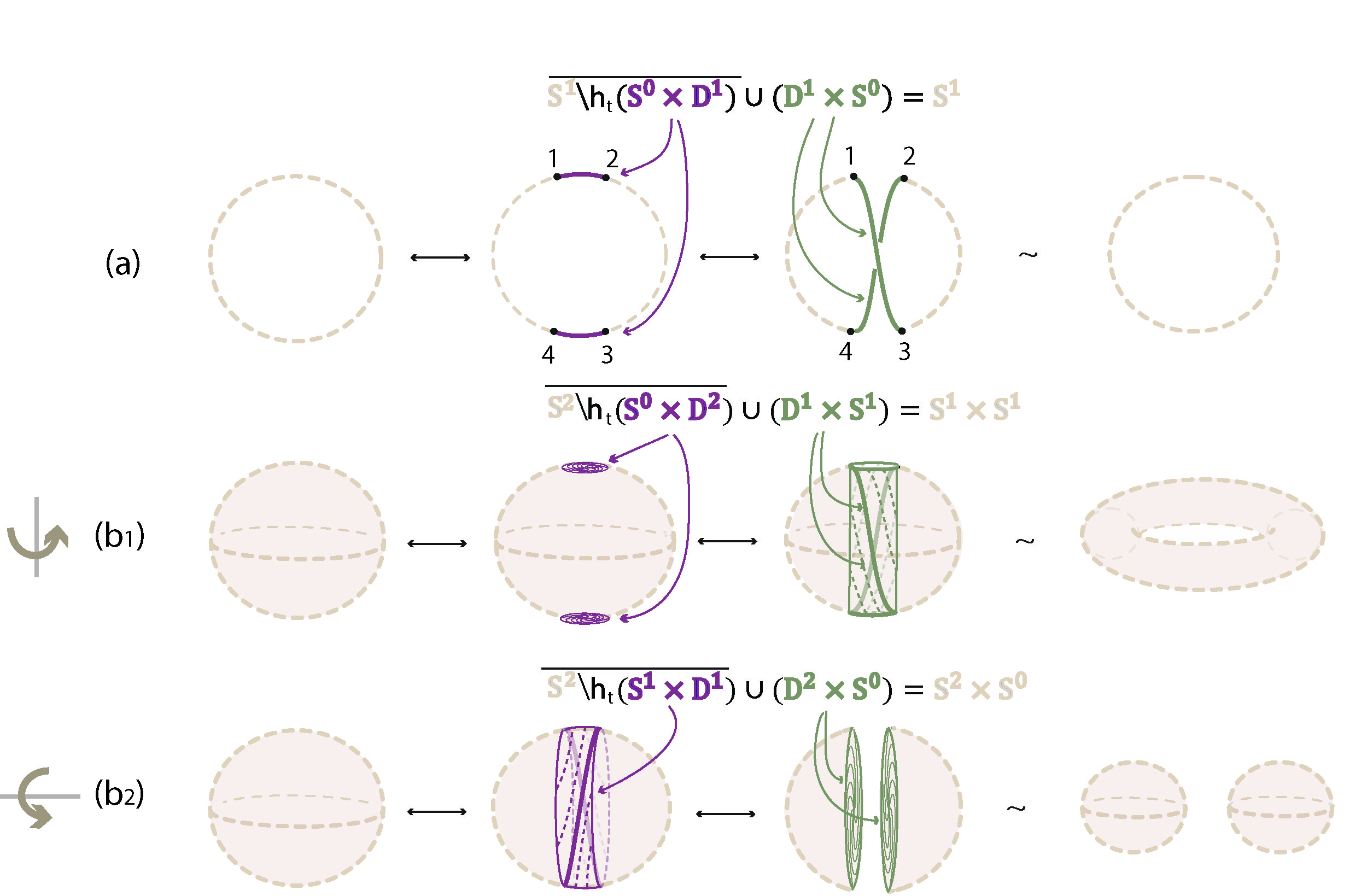}
\caption{Formal \textbf{(a)} 1-dimensional 0-surgery \textbf{(b\textsubscript{1})} 2-dimensional 0-surgery and \textbf{(b\textsubscript{2})} 2-dimensional 1-surgery using a twisting embedding $h_t$.}
\label{1_2D_Formal_twisted}\end{center}
\end{figure}

We only have one kind of surgery on a 1-manifold $M$,  the   \textit{1-dimensional 0-surgery}: \[M' = \chi(M) = \overline{M\setminus h(S^0\times D^{1})} \cup_{h|_{S^0\times S^{0}}} (D^{1}\times S^{0}). \] The above definition means that two segments $S^0\times D^1$ are removed from $M$ and they are replaced by two different segments $D^1 \times S^0$ by reconnecting the four boundary points $S^0\times S^0$ in a different way. In Fig.~\ref{1_2D_Formal} (a) and ~\ref{1_2D_Formal_twisted} (a),  $S^0\times S^0 =\{1, 2, 3, 4\}$. As one possibility,  if we start with $M=S^1$ and use as $h$ the standard (identity) embedding denoted with $h_s$,  we obtain two circles $S^1 \times S^0$. Namely,  denoting by $1$ the identity homeomorphism,  we have $h_s:S^0\times D^{1}=D^{1} \amalg D^{1} \xrightarrow{1 \amalg  1}  S^0\times D^{1} \hookrightarrow M$,  see Fig.~\ref{1_2D_Formal} (a). However,  we can also obtain one circle $S^1$ if $h$ is an embedding $h_t$ that reverses the orientation of one of the two arcs of $S^0\times D^1$. Then in the substitution,  joining endpoints 1 to 3 and 2 to 4,  the two new arcs undergo a half-twist,  see Fig.~\ref{1_2D_Formal_twisted} (a). More specifically,  if we take ${D^1}=[-1, +1]$ and define the homeomorphism  $\omega:D^{1}\to D^{1} ; t \to -t$,  the embedding used in Fig.~\ref{1_2D_Formal_twisted} (a) is $h_t:S^0\times D^{1}=D^{1} \amalg D^{1} \xrightarrow{1 \amalg  \omega}  S^0\times D^{1} \hookrightarrow M $ which rotates one $D^{1}$ by 180\si{\degree}. The difference between the embeddings $h_s$ and $h_t$ of $S^0\times D^1$ can be clearly seen by comparing the four boundary points $1, 2, 3$ and $4$ in Fig.~\ref{1_2D_Formal} (a) and Fig.~\ref{1_2D_Formal_twisted} (a).


\smallbreak

\smallbreak

Note that in dimension one,  the dual case is also an 1-dimensional 0-surgery. For example,  looking at the reverse process of Fig.~\ref{1_2D_Formal} (a),  we start with two circles $M'=S^1\amalg S^1$ and,  if each  segment of $D^1 \times S^0$ is embedded in a different circle,  the result of the (dual) 1-dimensional 0-surgery is one circle: $ \chi^{-1}(M')=M=S^1$.

\subsection{2-dimensional 0-surgery}\label{2D_FormalE0}
 Starting with a 2-manifold $M$,  there are two types of surgery. One type is the \textit{ 2-dimensional 0-surgery},  whereby two discs $S^0\times D^2$ are removed from  $M$ and are replaced in the closure of the remaining manifold by a cylinder $D^1\times S^1$,  which gets attached via a homeomorphism along the common boundary $S^0\times S^1$ comprising two copies of $S^1$. The gluing homeomorphism of the common boundary may twist one or both copies of $S^1$. For $M=S^2$ the above operation changes its homeomorphism type from the 2-sphere to that of the torus. View Fig.~\ref{1_2D_Formal}~({b\textsubscript{1}}) for the standard embedding  $h_s$ and Fig.~\ref{1_2D_Formal_twisted}~({b\textsubscript{1}}) for a twisting embedding $h_t$. For example,  the  homeomorphism  $\mu:D^{2}\to D^{2} ; (t_1, t_2) \to (-t_1, -t_2)$ induces the 2-dimensional analogue $h_t$ of the embedding defined in the previous example,  namely: $h_t:S^0\times D^{2}=D^{2} \amalg D^{2} \xrightarrow{1 \amalg  \mu}  S^0\times D^{2} \hookrightarrow M $ which rotates one $D^{2}$ by 180\si{\degree}. When,  now,  the cylinder $D^1\times S^1$ is glued along the common boundary  $S^0\times S^1$,  the twisting of this boundary induces the twisting of the cylinder,  see Fig.~\ref{1_2D_Formal_twisted}~({b\textsubscript{1}}). 

\subsection{2-dimensional 1-surgery}\label{2D_FormalE1}
The other possibility of 2-dimensional surgery on $M$ is the \textit{2-dimensional 1-surgery}: here a cylinder (or annulus) $S^1 \times D^1$ is removed from $M$ and is replaced in the closure of the remaining manifold by two discs  $D^2 \times S^0$ attached along the common boundary $S^1 \times S^0$. For $M=S^2$ the result is two copies of  $S^2$,  see Fig.~\ref{1_2D_Formal}~({b\textsubscript{2}}) for the standard embedding  $h_s$. Unlike Fig.~\ref{1_2D_Formal}~({b\textsubscript{1}})  where the cylinder is illustrated vertically,  in Fig.~\ref{1_2D_Formal}~({b\textsubscript{2}}),  the cylinder is illustrated horizontally. This choice was made so that the instances of 1-dimensional surgery can be obtained by crossections of the instances of both types of 2-dimensional surgeries,  see further Remark~\ref{Rot}. Fig.~\ref{1_2D_Formal_twisted}~({b\textsubscript{2}}) illustrates a twisting embedding $h_t$,  where a twisted cylinder is being removed. In that case,  taking $D^{1}=\{h:h \in [-1, 1]\}$ and homeomorphism $\zeta$: 

\begin{samepage} 
 \begin{center}
$\zeta:S^1\times D^{1}\to S^1\times D^{1};$ 
 \nopagebreak 
 \\[5pt] 
$\zeta: (t_1, t_2, h) \to (t_1\cos{\frac{(1-h)\pi}{2}}-t_2\sin{\frac{(1-h)\pi}{2}}, t_1\sin{\frac{(1-h)\pi}{2}}+t_2\cos{\frac{(1-h)\pi}{2}}, h)$
\end{center}  
 \end{samepage} 

{\noindent}the embedding $h_t$ is defined as: $h_t:S^1\times D^{1} \xrightarrow{\zeta} S^1\times D^{1} \hookrightarrow M $. This operation corresponds to fixing the circle $S^1$ bounding the right side of the cylinder $S^1 \times D^1$,  rotating the circle $S^1$ bounding the left side of the cylinder by 180\si{\degree} and letting the rotation propagate from left to right. This twisting of the cylinder can be seen by comparing the second instance of Fig.~\ref{1_2D_Formal}~({b\textsubscript{2}}) with  the second instance of Fig.~\ref{1_2D_Formal_twisted}~({b\textsubscript{2}}),  but also by comparing the third instance of Fig.~\ref{1_2D_Formal}~({b\textsubscript{1}}) with  the third instance of Fig.~\ref{1_2D_Formal_twisted}~({b\textsubscript{1}}).

It follows from Definition~\ref{surgery} that a dual 2-dimensional 0-surgery is a 2-dimensional 1-surgery and vice versa. Hence,  Fig.~\ref{1_2D_Formal}~({b\textsubscript{1}}) shows that a 2-dimensional 0-surgery on a sphere is the reverse process of a 2-dimensional 1-surgery on a torus. Similarly,  as illustrated in Fig.~\ref{1_2D_Formal}~({b\textsubscript{2}}),  a 2-dimensional 1-surgery on a sphere is the reverse process of a 2-dimensional 0-surgery on two spheres. In the figure the symbol $\longleftrightarrow $ depicts surgeries from left to right and their corresponding dual surgeries from right to left.


\begin{remark} \label{Rot} \rm  The stages of the process of 2-dimensional 0-surgery on $S^2$ can be obtained by rotating the stages of 1-dimensional 0-surgeries on $S^1$ by 180\si{\degree} around a vertical axis,  see Fig.~\ref{1_2D_Formal}~({b\textsubscript{1}}). Similarly,  the stages of 2-dimensional 1-surgery on $S^2$ can be obtained by rotating the stages of 1-dimensional 0-surgeries on $S^1$ by 180\si{\degree} around a horizontal axis,  see Fig.~\ref{1_2D_Formal}~({b\textsubscript{2}}). It follows from the above that 1-dimensional 0-surgery can be obtained as a cross-section of either type of 2-dimensional surgery.
\end{remark} 

\section{Dynamic topological surgery in natural processes}\label{Dynamics}
In this section we present natural processes exhibiting topological surgery in dimensions 1 and 2 and we incorporate the observed dynamics to a schematic model showing the intermediate steps that are missing from the formal definition. This model extends surgery to a continuous process caused by local forces. Note that these intermediate steps can also be explained through Morse theory,  see Remark~\ref{Morse} for details.

\subsection{Dynamic 1-dimensional topological surgery}\label{1D}
\begin{figure}[!h]
\begin{center}
\includegraphics[width=12cm]{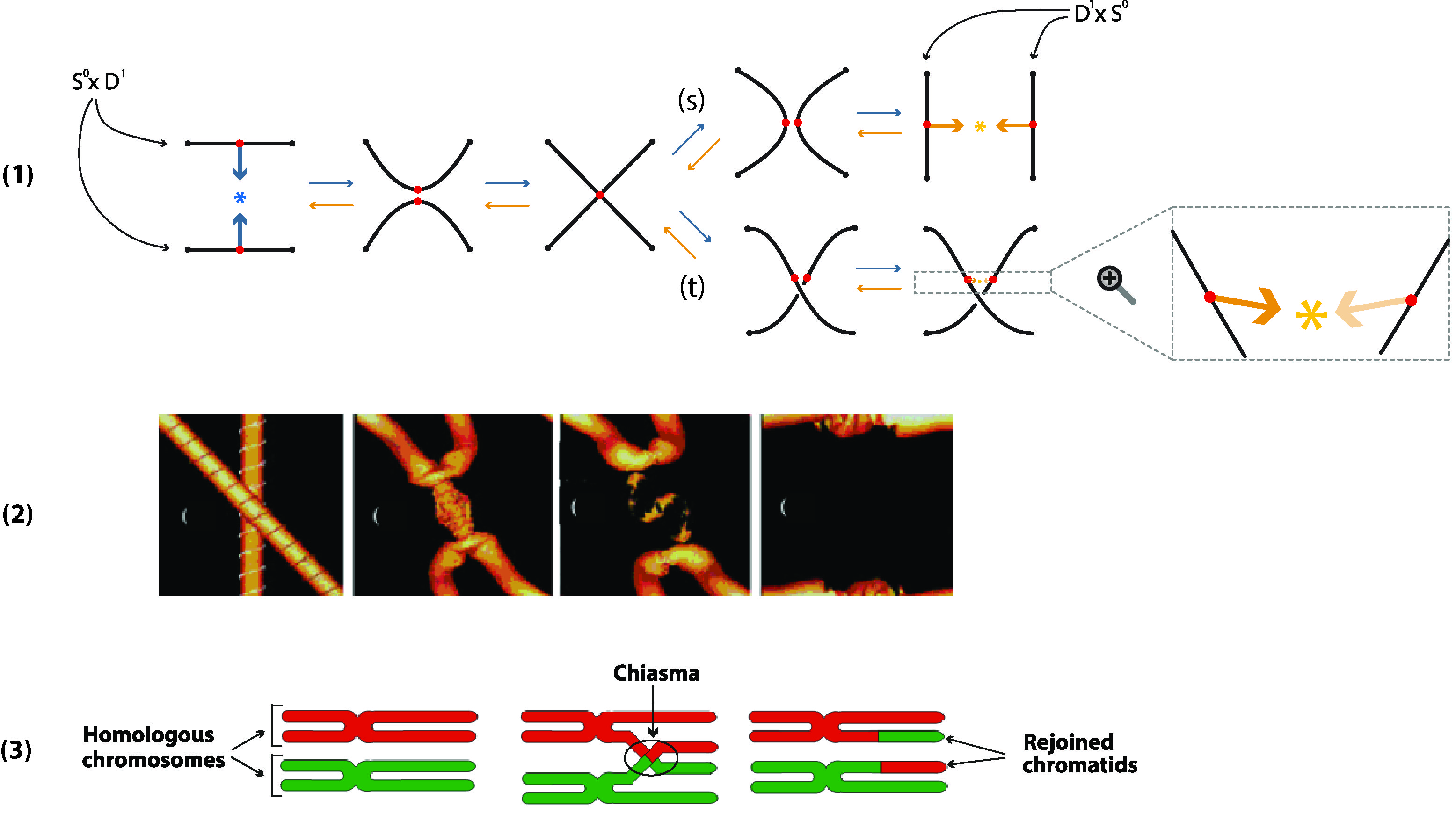}
\caption{\textbf{(1)} Dynamic 1-dimensional surgery locally \textbf{(2)} The reconnection of cosmic magnetic lines \textbf{(3)} Crossing over of chromosomes during meiosis.}
\label{1D_Nature}
\end{center}
\end{figure}

We find that $1$-dimensional 0-surgery is present in phenomena where 1-dimensional splicing and reconnection occurs. It can be seen for example during meiosis when new combinations of genes are produced,  see Fig.~\ref{1D_Nature} (3),   and during magnetic reconnection,  the phenomena whereby cosmic magnetic field lines from different magnetic domains are spliced to one another,  changing their pattern of conductivity with respect to the sources,  see Fig.~\ref{1D_Nature} (2) from \cite{DaAn}. It is worth mentioning that $1$-dimensional 0-surgery is also present during the reconnection of  vortex tubes in a viscous fluid and quantized vortex tubes in superfluid helium. As mentioned in  \cite{LaRiSu},  these cases have some common qualitative features with the magnetic reconnection shown in Fig.~\ref{1D_Nature} (2).

In fact,  all the above phenomena have similar dynamics. Namely,  $1$-dimensional 0-surgery is a  continuous process for all of them. Furthermore,  as detailed in \cite{SSN1},  in most cases,  surgery is the result of local forces. These common features are captured by our model in Fig.~\ref{1D_Nature} (1) which describes the process of dynamic $1$-dimensional 0-surgery locally. The process starts with the two points (in red) specified on any 1-dimensional manifold,  on which attracting forces are applied (in blue). We assume that these forces are caused by an attracting center (also in blue). Then,  the two segments $S^0\times D^1$,  which are neighborhoods of the two points,  get close to one another. When the specified points (or centers) of the two segments reach the attracting center they touch and recoupling takes place,  giving rise to the two final segments $D^1 \times S^0$,  which split apart. In Fig.~\ref{1D_Nature} (1),  case (s) corresponds to the identity embedding $h_s$ described in Section~\ref{1D_FormalE},  while (t) corresponds to the  twisting embedding $h_t$ described in Section~\ref{1D_FormalE}.

\smallbreak
 
As also mentioned in Section~\ref{1D_FormalE},  the dual case is also a 1-dimensional 0-surgery,  as it removes segments $D^1 \times S^0$ and replaces them by segments $S^0\times D^1$. This is the reverse process which starts from the end and is illustrated in Fig.~\ref{1D_Nature} (1) as a result of the orange forces and attracting center which are applied on the `complementary' points.

Finally,  for details about the described dynamics and attracting forces in the aforementioned phenomena,  the reader is referred to the analysis done in \cite{SSN1}.

\begin{remark}\label{Morse} \rm
It is worth mentioning that the intermediate steps of surgery presented in Fig.~\ref{1D_Nature} (1) can also be viewed in the context of Morse theory \cite{Mil}. By using the local form of a Morse function,  we can visualize the process of surgery by varying parameter $t$ of equation $x^2-y^2=t$. For $t=-1$ it is the hyperbola shown in the second stage of Fig.~\ref{1D_Nature} (1) where the two segments get close to one another. For $t=0$ it is the two straight lines where the reconnection takes place as shown in the third stage of Fig.~\ref{1D_Nature} (1) while for $t=1$ it represents the hyperbola of the two final segments shown in case (s) of the fourth stage of Fig.~\ref{1D_Nature} (1). This sequence can be generalized for higher dimensional surgeries as well,  however,  in this paper we will not use this approach as we are focusing on the introduction of forces and of the attracting center. 
\end{remark}

\subsection{Dynamic 2-dimensional topological surgery}\label{2D}
\begin{figure}[!h]
\begin{center}
\includegraphics[width=12cm]{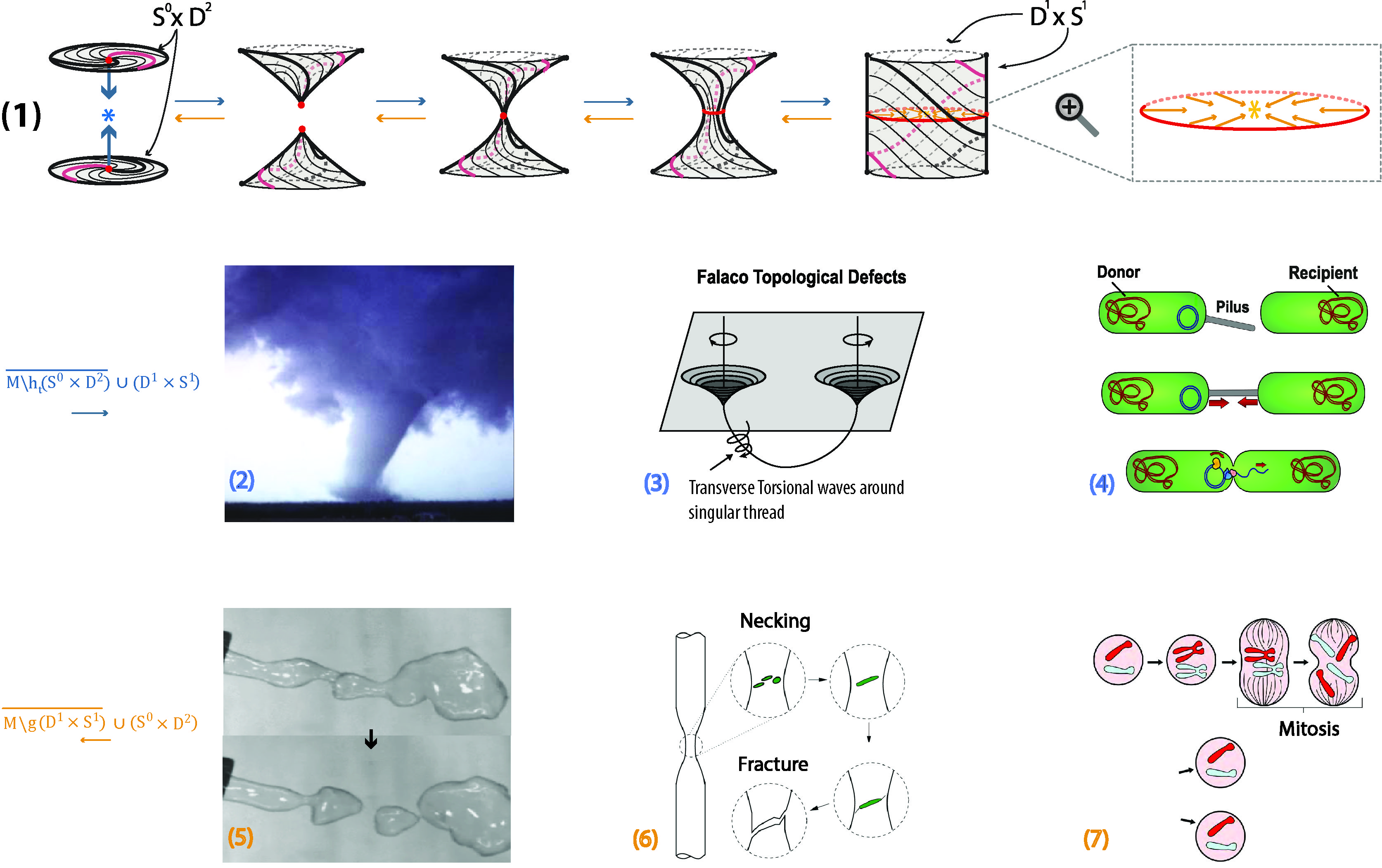}
\caption{\textbf{(1)} Dynamic 2-dimensional surgery locally \textbf{(2)} Tornadoes \textbf{(3)} Falaco solitons  \textbf{(4)} Gene transfer in bacteria \textbf{(5)} Soap bubble splitting \textbf{(6)} Fracture \textbf{(7)} Mitosis}
\label{2D_Nature}
\end{center}
\end{figure}

We find that $2$-dimensional surgery is present in phenomena where 2-dimensional merging and recoupling occurs. For example,  \textit{2-dimensional 0-surgery} can be seen during the formation of tornadoes,   see Fig.~\ref{2D_Nature} (2) (this phenomenon will be detailed in Section~\ref{E2D0_D3S0}). Further,  it can be seen in the formation of Falaco solitons,  see Fig.~\ref{2D_Nature} (3) (note that each Falaco soliton consists of a pair of locally unstable but globally stabilized contra-rotating identations in the water-air discontinuity surface of a swimming pool,  see \cite{Ki} for details).
It can also be seen in gene transfer in bacteria where the donor cell produces a connecting tube called a `pilus' which attaches to the recipient cell,  see Fig.~\ref{2D_Nature} (4); in drop coalescence,  the phenomenon where two dispersed drops merge into one,  and in the formation of black holes (this phenomena will be discussed in Section~\ref{E2D0}),  see Fig.~\ref{ES2D0_B3} (2) (Source: www.black-holes.org).  

On the other hand,  \textit{ 2-dimensional 1-surgery} can be seen during soap bubble splitting,  where a soap bubble splits into two smaller bubbles,  see Fig.~\ref{2D_Nature} (5) (Source: soapbubble.dk); when the tension applied on metal specimens by tensile forces results in the phenomena of necking and then fracture,  see Fig.~\ref{2D_Nature} (6); also in the biological process of mitosis,  where a cell splits into two new cells,  see Fig.~\ref{2D_Nature} (7).

\smallbreak

Phenomena exhibiting 2-dimensional 0-surgery are the results of \textit{two colinear attracting forces which `create' a cylinder}. These phenomena have similar dynamics and are characterized by their continuity and the attracting forces causing them,  see \cite{SSN1} for details. These common features are captured by our model in Fig.~\ref{2D_Nature} (1) which describes both cases of dynamic $2$-dimensional surgeries locally. Note that we have used non-trivial embeddings (recall Section~\ref{2D_FormalE0}) which are more appropriate for natural processes involving twisting,  such as tornadoes and Falaco solitons. In Fig.~\ref{2D_Nature} (1),  both initial discs have been rotated in opposite directions by an angle of $3\pi/4$ and we can see how this rotation induces the twisting of angle $3\pi/2$ of the final cylinder.


The process of \textit{dynamic $2$-dimensional 0-surgery} starts with two points,  or poles,  specified on the manifold (in red) on which attracting forces caused by an attracting center are applied (in blue).  Then,  the two discs $S^0\times D^2$,   neighbourhoods of the two poles,  approach each other. When the centers of the two discs touch,  recoupling takes place and the discs get transformed into the final cylinder $D^1\times S^1$,  see Fig.~\ref{2D_Nature} (1). The cylinder created during 2-dimensional 0-surgery can take various forms. For example,  it is a tubular vortex of air in the case of tornadoes,  a transverse torsional wave in the case of Falaco solitons and a pilus joining the genes in gene transfer in bacteria.

On the other hand,  phenomena exhibiting 2-dimensional 1-surgery are the result of \textit{an infinitum of coplanar attracting forces which `collapse' a cylinder},  see Fig.~\ref{2D_Nature} (1) from the end. As mentioned in Section~\ref{2D_FormalE1},  the dual case of 2-dimensional 0-surgery is the 2-dimensional 1-surgery and vice versa. This is illustrated in Fig.~\ref{2D_Nature} (1) where the reverse process is the \textit{2-dimensional 1-surgery} which starts with the cylinder and a specified circular region (in red) on which attracting forces caused by an attracting center are applied (in orange). A `necking'  occurs in the middle which degenerates into a point and finally tears apart creating two discs $S^0\times D^2$. This cylinder can be embedded,  for example,  in the region of the bubble's  surface where splitting occurs,  on the region of metal specimens where necking and fracture occurs or on the  equator of the cell which is about to undergo a mitotic process.

\begin{remark} \label{Rot2} \rm  From Definition~\ref{surgery} we know that surgery always starts by removing a thickened sphere $S^n\times D^{m-n}$ from the initial manifold. As seen in  Fig.~\ref{2D_Nature} (1),  in the case of 2-dimensional 0-surgery,  forces (in blue) are applied on two points,  or $S^0$,  whose thickening comprises the two discs,  while in the case of the 2-dimensional 1-surgery,  forces (in orange) are applied on a circle $S^1$,  whose thickening is the cylinder. In other words the forces that model a $2$-dimensional $n$-surgery are always applied to the core $n$-embedding $e=h_{|}: {S^n}={S^n} \times \{0\} \hookrightarrow  M$ of the framed $n$-embedding $h:S^n\times D^{2-n}\hookrightarrow  M$. Also,  note that Remark~\ref{Rot} is also true here. One can obtain Fig.~\ref{2D_Nature} (1) by rotating Fig.~\ref{1D_Nature} (1) and this extends also to the dynamics and forces. For instance,  by rotating the two points,  or $S^0$,  on which the pair of forces of 1-dimensional 0-surgery acts (shown in red in the last instance of Fig.~\ref{1D_Nature} (1)) by 180\si{\degree} around a vertical axis we get the circle,  or $S^1$,  on which the infinitum of coplanar attracting forces of 2-dimensional 1-surgery acts (shown in red in the last instance of Fig.~\ref{2D_Nature} (1)).
\end{remark} 

Finally,  it is worth pointing out that these local dynamics produce different manifolds depending on the initial manifold where they act. For example,  2-dimensional 0-surgery transforms an $S^0\times S^2$ to an $S^2$ by adding a cylinder during gene transfer in bacteria (see Fig.~\ref{2D_Nature} (4)) but can also transform an $S^2$ to a torus by `drilling out' a cylinder during the formation of Falaco solitons (see Fig.~\ref{2D_Nature} (3)) in which case $S^2$ is the pool of water and the cylinder is the boundary of the tubular neighborhood around the thread joining the two poles. 


\section{Defining solid topological surgery} \label{DSolid}
Looking closer at the phenomena exhibiting 2-dimensional surgery shown in Fig.~\ref{2D_Nature},  one can see that,  with the exception of soap bubble splitting that involves surfaces,  all others involve 3-dimensional manifolds. For instance,  what really happens during a mitotic process is that a solid cylindrical region located in the center of the cell collapses and a $D^3$ is transformed into an $S^0\times D^3$. Similarly,  during tornado formation,  the created cylinder is not just a cylindrical surface  $D^1\times S^1$ but a solid cylinder $D^2\times S^1$ containing many layers of air (this phenomena will be detailed in Section~\ref{E2D0_D3S0} ). Of course we can say that,  for phenomena involving 3-dimensional manifolds,  the outer layer of the initial manifold is undergoing 2-dimensional surgery. In this section we will define topologically what happens to the whole manifold. 

The need of such a definition is also present in dimension 1 for modeling phenomena such as the merging of oil slicks and tension on membranes (or soap films). These phenomena undergo the process of 1-dimensional 0-surgery but happen on surfaces instead of 1-manifolds. 

We will now introduce the notion of solid surgery (in both dimensions 1 and 2) where the interior of the initial manifold is filled in. There is one key difference compared to the dynamic surgeries discussed in the previous section. While the  local dynamics described in Fig.~\ref{1D_Nature} and ~\ref{2D_Nature} can be embedded in any manifold,  here we also have to fix the initial manifold in order to define solid surgery. For example,  as we will see next,  we define separately the processes of \textit{solid 1-dimensional 0-surgery on $D^2$} and \textit{ solid 1-dimensional 0-surgery on $D^2 \times S^0$}. However,  the underlying features are common in both.

\subsection{Solid 1-dimensional topological surgery}\label{Solid1D}
\begin{figure}[!h]
\begin{center}
\includegraphics[width=12cm]{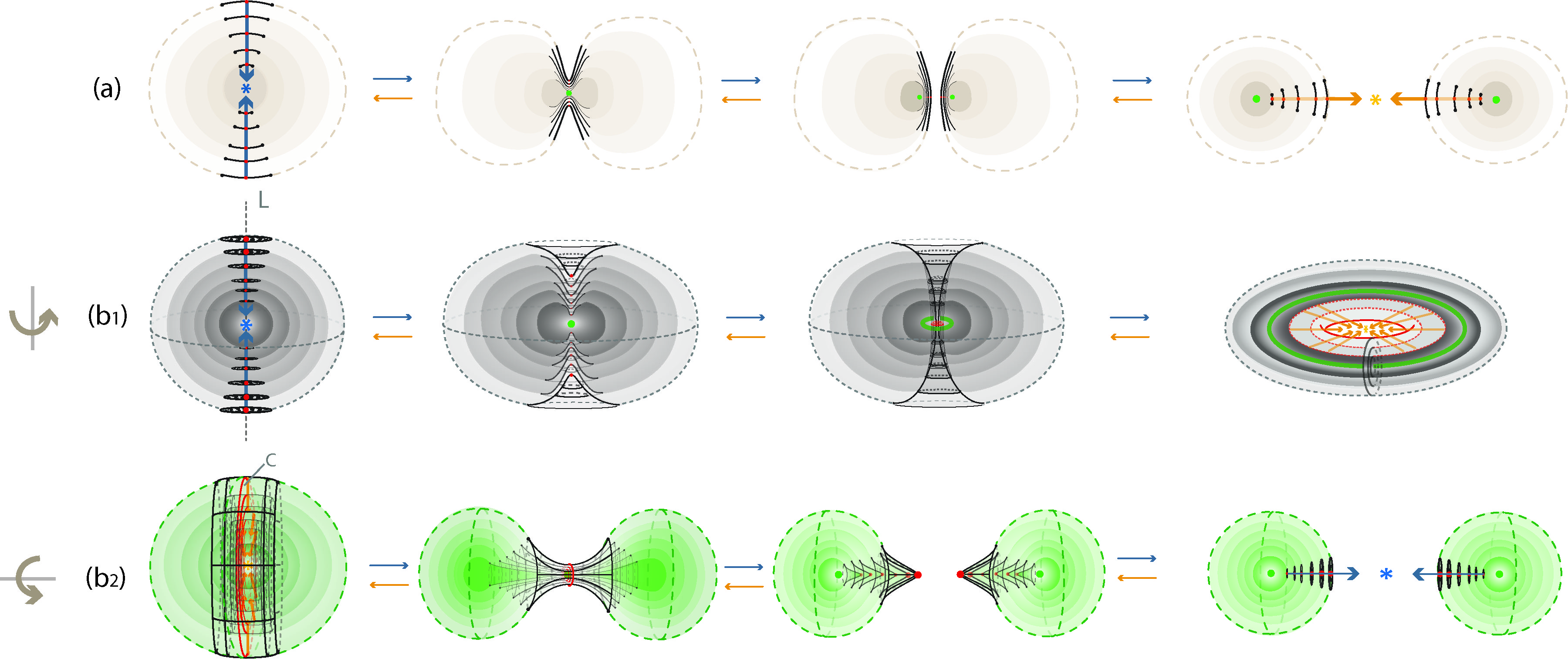}
\caption{\textbf{(a)} Solid 1-dimensional 0-surgery on $D^2$ \textbf{(b\textsubscript{1})} Solid 2-dimensional 0-surgery on on $D^3$  \textbf{(b\textsubscript{2})} Solid 2-dimensional 1-surgery on $D^3$}
\label{Solid}
\end{center}
\end{figure}

The process of solid 1-dimensional 0-surgery on $D^2$ is equivalent to performing 1-dimensional 0-surgeries on the whole continuum of concentric circles included in $D^2$,  see Fig.~\ref{Solid} (a). More precisely,  and introducing at the same time dynamics,  we define:

\smallbreak

\begin{definition} \rm  \textit{ Solid 1-dimensional 0-surgery on $D^2$} is the following process. We start with the $2$-disc of radius 1 with polar layering: 
$$
D^2 = \cup_{0<r\leq 1} S^1_r \cup \{P\}, 
$$ 
where $r$ the radius of a circle and $P$ the limit point of the circles,  which is the center of the disc and also the circle of radius zero. We specify colinear pairs of antipodal points,  all on the same diameter,  with neighbourhoods of analogous lengths,  on which the same colinear attracting forces act. See Fig.~\ref{Solid} (a) where these forces and the attracting center are shown in blue. Then the antipodal segments get closer to one another or,  equivalently,  closer to the attracting center. Note that here,  the attracting center coincides with the limit point of all concentric circles,  which is shown in green from the second instance and on. Then,  we perform 1-dimensional 0-surgery on the whole continuum of concentric circles.  We define 1-dimensional 0-surgery on the limit point $P$ to be the two limit points of the resulting surgeries. That is,  the effect of \textit{solid 1-dimensional 0-surgery on a point is the creation of two new points}. Next,  the segments reconnect until we have two copies of $D^2$. 
 
The above process is the same as first removing the center $P$ from $D^2$,  doing the  1-dimensional 0-surgeries and then taking the closure of the resulting space. The resulting manifold is 
$$
\chi(D^2) := \cup_{0<r\leq 1}\chi(S^1_r) \cup \chi(P), 
$$
which comprises two copies of $D^2$.

We also have the reverse process of the above,  namely,   \textit{ solid 1-dimensional 0-surgery on two discs $D^2 \times S^0$}. This process is the result of the orange forces and attracting center which are applied on the `complementary' points,  see Fig.~\ref{Solid} (a) in reverse order. This operation is equivalent to performing 1-dimensional 0-surgery on the  whole continuum of pairs of concentric circles in $D^2\amalg D^2$. We only need to define solid 1-dimensional 0-surgery on two limit points to be the limit point $P$ of the resulting surgeries. That is,  the effect of \textit{ solid 1-dimensional 0-surgery on two points is their merging into one point}. The above process is the same as first removing the centers from the $D^2 \times S^0$,  doing the 1-dimensional 0-surgeries and then taking the closure of the resulting space. The resulting manifold is 
$$
\chi^{-1}(D^2 \times S^0) := \cup_{0<r\leq 1}\chi^{-1}(S^1_r \times S^0) \cup \chi^{-1}(P \times S^0), 
$$
which comprises one copy of $D^2$.
\end{definition}

\subsection{Solid 2-dimensional topological surgery}\label{Solid2D}
Moving up one dimension,  there are two types of solid 2-dimensional surgery on the $3$-ball,  $D^3$,  analogous to the two types of 2-dimensional surgery. 
More precisely we have:

\begin{definition} \label{cont2D}  \rm We start with the $3$-ball of radius 1 with polar layering: 
$$
D^3 = \cup_{0<r\leq 1} S^2_r \cup \{P\}, 
$$ 
where $r$ the radius of the 2-sphere $S^2_r$ and $P$ the limit point of the spheres,  that is,  their common center and the center of the ball. \textit{Solid 2-dimensional 0-surgery on $D^3$} is the  topological procedure whereby 2-dimensional 0-surgery takes place on each spherical layer that $D^3$ is made of. More precisely,  as illustrated in Fig.~\ref{Solid}~({b\textsubscript{1}}),  on all spheres $S^2_r$  colinear pairs of antipodal points are specified,  all on the same diameter,  on which the same colinear attracting forces act. The poles have disc neighborhoods of analogous areas. Then,  2-dimensional 0-surgeries are performed on the whole continuum of the concentric spheres using the same  embedding $h$ (recall Section~\ref{2D_FormalE0}). Moreover,  2-dimensional 0-surgery on the limit point $P$ is defined to be the limit circle of the nested tori resulting from the continuum of 2-dimensional 0-surgeries. That is,  the effect of \textit{2-dimensional 0-surgery on a point is defined to be the creation of a circle}.

The process is characterized on one hand by the 1-dimensional core $L$ of the solid cylinder which joins the two selected antipodal points of the outer shell and intersects each spherical layer at its two corresponding antipodal points,  and on the other hand by the embedding $h$. The process results in a continuum of layered tori and can be viewed as drilling out a tunnel along $L$ according to $h$. Note that in Fig.~\ref{Solid},  the identity embedding has been used. However,  a twisting embedding,  which is the case shown in Fig.~\ref{2D_Nature} (1),  agrees with our intuition that,  for opening a hole,  \textit{drilling with twisting} seems to be the easiest way. Examples of these two embeddings can be found in Section~\ref{2D_FormalE0}.

Furthermore,  \textit{solid 2-dimensional 1-surgery on $D^3$} is the  topological procedure where on all spheres $S^2_r$ nested cylindrical peels of the solid cylinder of analogous areas are specified and the same coplanar attracting forces act on all spheres,  see  Fig.~\ref{Solid}~({b\textsubscript{2}}). Then,  2-dimensional 1-surgeries are performed on the whole continuum of the concentric spheres using the same embedding $h$. Moreover,  2-dimensional 1-surgery on the limit point $P$  is defined to be the two limit points of the nested pairs of 2-spheres resulting from  the continuum of 2-dimensional surgeries. That is,  the effect of \textit{2-dimensional 1-surgery on a point is the creation of two new points}.
The process is characterized by the 2-dimensional central disc of the solid cylinder and the embedding $h$,  and it can be viewed as squeezing the central disc $C$ or,  equivalently,  as pulling apart the left and right hemispheres with possible twists,  if $h$ is a twisting embedding. This agrees with our intuition that for cutting a solid object apart,  \textit{pulling with twisting} seems to be the easiest way. Examples of the identity and the twisting embedding can be found in Section~\ref{2D_FormalE1}. 

For both types of solid 2-dimensional surgery,  the above process is the same as: first removing the center $P$ from $D^3$,  performing the 2-dimensional surgeries and then taking the closure of the resulting space. Namely we obtain:  
$$
\chi(D^3) := \cup_{0<r\leq 1}\chi(S^2_r) \cup \chi(P), 
$$
which is a solid torus in the case of solid 2-dimensional 0-surgery and  two copies of $D^3$ in the case of solid 2-dimensional 1-surgery.
\end{definition}

\smallbreak
As seen in Fig.~\ref{Solid},  we also have the two dual solid 2-dimensional surgeries,  which represent the reverse processes. As already mentioned in Section~\ref{2D_FormalE1},  the dual case of 2-dimensional 0-surgery is the 2-dimensional 1-surgery and vice versa. More precisely:

\begin{definition} \label{solto2D} \rm The \textit{dual case of solid 2-dimensional 0-surgery on $D^3$} is the solid 2-dimensional 1-surgery on a solid torus $D^2 \times S^1$. This is the reverse process shown in Fig.~\ref{Solid}~({b\textsubscript{1}}) which results from the orange forces and attracting center. Given that the solid torus can be written as a union of nested tori together with the core circle: $D^2 \times S^1=(\cup_{0<r\leq 1}S^1_r \cup \{0\}) \times S^1$,  solid 2-dimensional 1-surgeries are performed on each toroidal layer starting from specified annular peels of analogous sizes where the same coplanar forces act on the central rings of the annuli. These forces are caused by the same attracting center lying outside the torus. It only remain to define the solid 2-dimensional 1-surgery on the limit circle to be the limit point $P$ of the resulting surgeries. That is,  the effect of \textit{ solid 2-dimensional 1-surgery on the core circle is that it collapses into one point,  the attracting center}. The above process is the same as first removing the core circle from $D^2 \times S^1$,  doing the 2-dimensional 1-surgeries on the layered tori,  with the same coplanar acting forces,  and then taking the closure of the resulting space. Hence,  the resulting manifold is 

$$
\chi^{-1}(D^2 \times S^1) := \cup_{0<r\leq 1}\chi^{-1}(S^1_r \times S^1) \cup \chi^{-1}(\{0\} \times S^1), 
$$
which comprises one copy of $D^3$.

Further,  \textit{the dual case of solid 2-dimensional 1-surgery on $D^3$} is the solid 2-dimensional 0-surgery on two $3$-balls $D^3$. This is the reverse process shown in Fig.~\ref{Solid}~({b\textsubscript{2}}) which results from the blue forces and attracting center. We only need to define the solid 2-dimensional 0-surgery on two limit points to be the limit point $P$ of the resulting surgeries. That is,  as in solid 1-dimensional surgery (see Fig.~\ref{Solid} (a)),  the effect of \textit{ solid 2-dimensional 0-surgery on two points is their merging into one point}. The above process is the same as first removing the centers from the $D^3 \times S^0$,  doing the 2-dimensional 0-surgeries on the nested spheres with the same colinear forces and then taking the closure of the resulting space. The resulting manifold is 
$$
\chi^{-1}(D^3 \times S^0) := \cup_{0<r\leq 1}\chi^{-1}(S^2_r \times S^0) \cup \chi^{-1}(P \times S^0), 
$$
which comprises one copy of $D^3$.
\end{definition}

Note that remarks~\ref{Rot} and~\ref{Rot2} are also true here. One can obtain Fig.~\ref{Solid} (b\textsubscript{1}) and  Fig.~\ref{Solid} (b\textsubscript{2}) by rotating Fig.~\ref{Solid} (a) respectively by 180\si{\degree} around a vertical axis and by 180\si{\degree} around a horizontal axis.

\begin{remark} \rm The notions of 2-dimensional (resp. solid 2-dimensional) surgery,  can be generalized from $S^2$ (resp. $D^3$) to a surface (resp. a handlebody) of genus $g$ creating a surface (resp. a handlebody) of genus $g \pm 1$ or a disconnected surface (resp. handlebody).
\end{remark}

\subsection{A dynamical system modeling solid 2-dimensional 0-surgery}\label{DS}

In \cite{SSN1, SSN2} we present and analyze a dynamical system for modeling 2-dimensional 0-surgery. This system was introduced in \cite{SaGr1} and is a generalization of the classical Lotka--Volterra system in three dimensions where the classic predator-prey population model is generalized to a two-predator and one-prey system. The parameters of the system affect the dynamics of the populations and they are analyzed in order to determine the bifurcation properties of the system.  

\smallbreak
\begin{figure}[!h]
\begin{center}
\includegraphics[width=12cm]{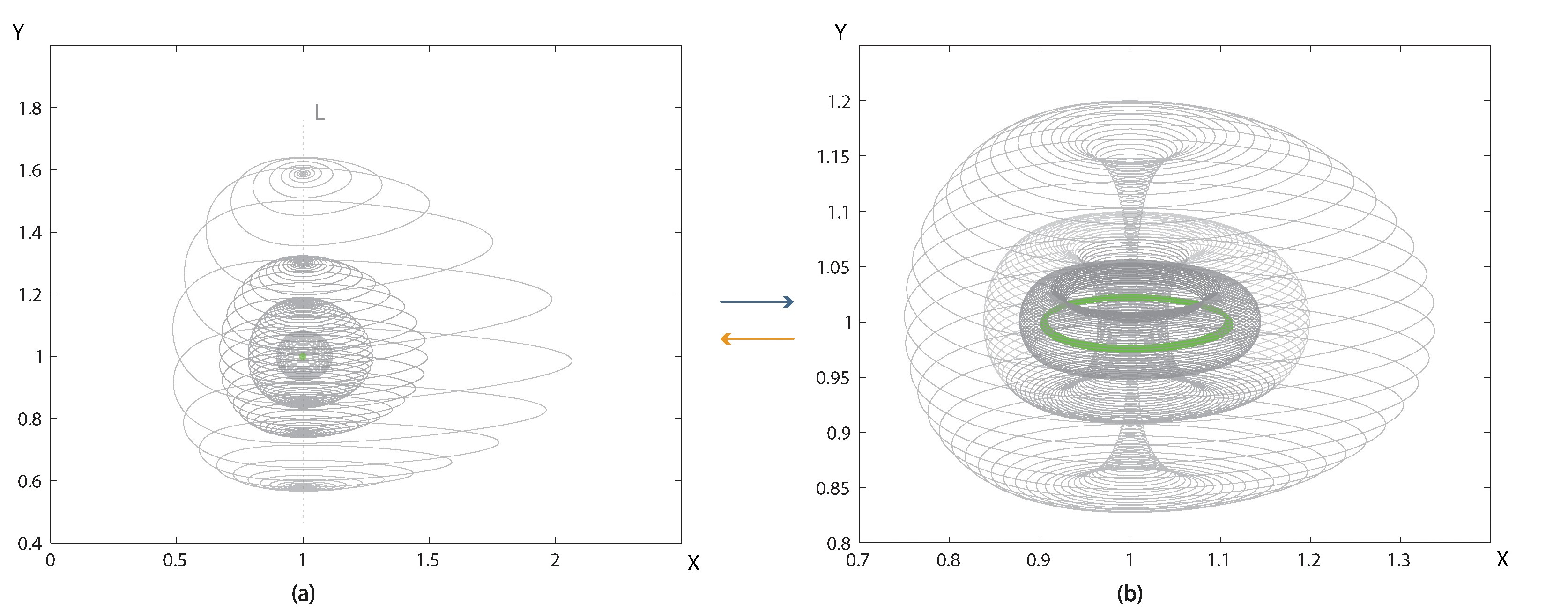}
\caption{Solid 2-dimensional 0-surgery by changing parameter space from (a) to (b).}
\label{LV_Surgery}
\end{center}
\end{figure}
\smallbreak 

In Fig.~\ref{LV_Surgery} we reproduce the numerical simulations done in \cite{SaGr1, SaGr2},  illustrating how by  changing the parameters space,  solid 2-dimensional 0-surgery is performed on the trajectories of the system. We have used the same colors as in Fig.~\ref{Solid}~({b\textsubscript{1}}) to make comparison easier. What is even more striking is that the new topological definition stating that the effect of 2-dimensional 0-surgery on a point is defined to be the creation of a circle (recall Section~\ref{Solid2D}) also has a meaning in the language of dynamical systems. More precisely,  the limit point in the spherical nesting of trajectories shown in green in Fig.~\ref{LV_Surgery} (a) is a steady state point and the core of the toroidal nesting of trajectories shown in green in Fig.~\ref{LV_Surgery} (b) is a limit cycle. This connection is detailed in \cite{SSN1},  see also \cite{SSN2}.

\section{Defining embedded topological surgery} \label{Embedded}

As mentioned in the Introduction,  we noticed that the ambient space is also involved in some natural processes exhibiting surgery. As we will see in this section,  depending on the dimension of the manifold,  the ambient space either leaves `room' for the initial manifold to assume a more complicated configuration or it participates more actively in the process. Independently of dimensions,  embedding surgery has the advantage that it allows us to view surgery as a process happening inside a space instead of abstractly. We define it as follows:  

\begin{definition} \rm   \textit{An embedded $m$-dimensional $n$-surgery} is a $m$-dimensional $n$-surgery where the initial manifold is an $m$-embedding $e:M \hookrightarrow S^d$,  $d\geq m$ of some $m$-manifold $M$,  and the result is also viewed as embedded in $S^d$. Namely,  according to Definition~\ref{surgery}: 
\[M' = \chi(e(M)) = \overline{e(M)\setminus h(S^n\times D^{m-n})} \cup_{h|_{S^n\times S^{m-n-1}}} D^{n+1}\times S^{m-n-1}  \hookrightarrow S^d.  \]
\end{definition}

{\noindent}Since in this analysis we focus on phenomena exhibiting embedded 1- and 2-dimensional surgery in 3-space,  from now on we fix $d=3$ and,  for our purposes,  we consider $S^3$ or ${\mathbb R}^3$ as our standard 3-space.

\subsection{Embedded 1-dimensional surgery} \label{Embedded1D}

In dimension 1,  the notion of embedded surgery allows the topological modeling of phenomena with more complicated initial 1-manifolds. Let us demonstrate this with the example of site-specific {\bf DNA recombination}. In this process,  the initial manifold is a (circular or linear) DNA molecule. With the help of certain enzymes,  site-specific recombination performs a 1-dimensional 0-surgery on the molecule,  causing possible knotting or linking of the molecule.

\smallbreak
\begin{figure}[!h]
\sidecaption[t]
\includegraphics[width=7cm]{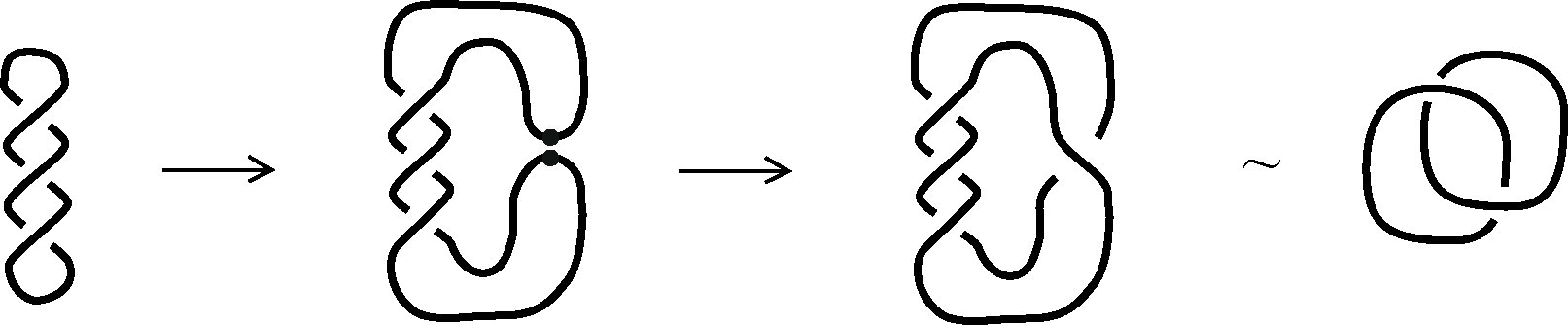}
\caption{DNA Recombination as an example of embedded 1-dimensional 0-surgery.}
\label{1D_Embedded}
\end{figure}

The first electron microscope picture of knotted DNA was presented in \cite{WaDuCo}. In this experimental study,  we see how genetically engineered circular DNA molecules can form DNA knots and links through the action of a certain recombination enzyme. A similar picture is presented in Fig.~\ref{1D_Embedded},  where site-specific recombination of a DNA molecule produces the Hopf link. It is worth mentioning that there are infinitely many knot types and that 1-dimensional 0-surgery on a knot may change the knot type or even result in a two-component link (as shown in Fig.~\ref{1D_Embedded}).  Since a knot is by definition an embedding of $M=S^1$ in $S^3$ or ${\mathbb R}^3$,  in this case embedded 1-dimensional surgery is the so-called \textit {knot surgery}. A good introductory book on knot theory is \cite{Ad} among many others.

We can summarize the above by stating that for $M=S^1$,  embedding in $S^3$ allows the initial manifold to become any type of knot. More generally,  in dimension 1 the ambient space which is of codimension 2 gives enough `room' for the initial 1-manifold to assume a more complicated homeomorphic configuration. 

\begin{remark}\label{EmbededSolid1D} \rm
Of course we also have,  in theory,  the notion of embedded solid 1-dimensional 0-surgery whereby the initial manifold is an embedding of a disc in 3-space. 


\end{remark}

\subsection{Embedded 2-dimensional surgery} \label{Embedded2D}

Passing now to 2-dimensional surgeries,  let us first note that an embedding of a sphere $M=S^2$ in $S^3$ presents no knotting because knotting requires embeddings of codimension 2. However,  in this case the ambient space plays a different role. Namely,  embedding 2-dimensional surgeries {\it allows the complementary space of the initial manifold to participate actively in the process}. Indeed,  while some natural phenomena undergoing surgery can be viewed as `local',  in the sense that they can be considered independently from the surrounding space,  some others are intrinsically related to the surrounding space. This relation can be both \textit{causal},  in the sense that the ambient space is involved in the triggering of the forces causing surgery,  and \textit{consequential},  in the sense that the forces causing surgery,  can have an impact on the ambient space in which they take place. 

Let us recall here that the ambient space $S^3$ can be viewed as ${\mathbb R}^3$ with all points at infinity compactified to one single point: $S^3 = {\mathbb R}^3 \cup \{\infty\}$. Further,  it can be viewed as the union of two $3$-balls along the common boundary: $S^3 = B^3 \cup D^3$ where a neighbourhood of the point at infinity can stand for one of the two $3$-balls. Finally,  $S^3$ can be viewed as the union of two solid tori along their common boundary: $S^3 = V_1\cup V_2$.

As mentioned in the introduction of Section~\ref{DSolid},  in most natural phenomena that exhibit 2-dimensional surgery,  the initial manifold is a \textit{solid} 3-dimensional object. Hence,  in the next subsections,  we describe natural phenomena undergoing solid 2-dimensional surgeries which exhibit the causal or consequential relation to the ambient space mentioned above and are therefore better described by considering them as embedded in $S^3$ or in ${\mathbb R}^3$. In parallel,  we describe how these processes are altering the whole space $S^3$ or ${\mathbb R}^3$.

\subsubsection{A topological model for black hole formation} \label{E2D0}

\begin{figure}[!h]
\begin{center}
\includegraphics[width=11cm]{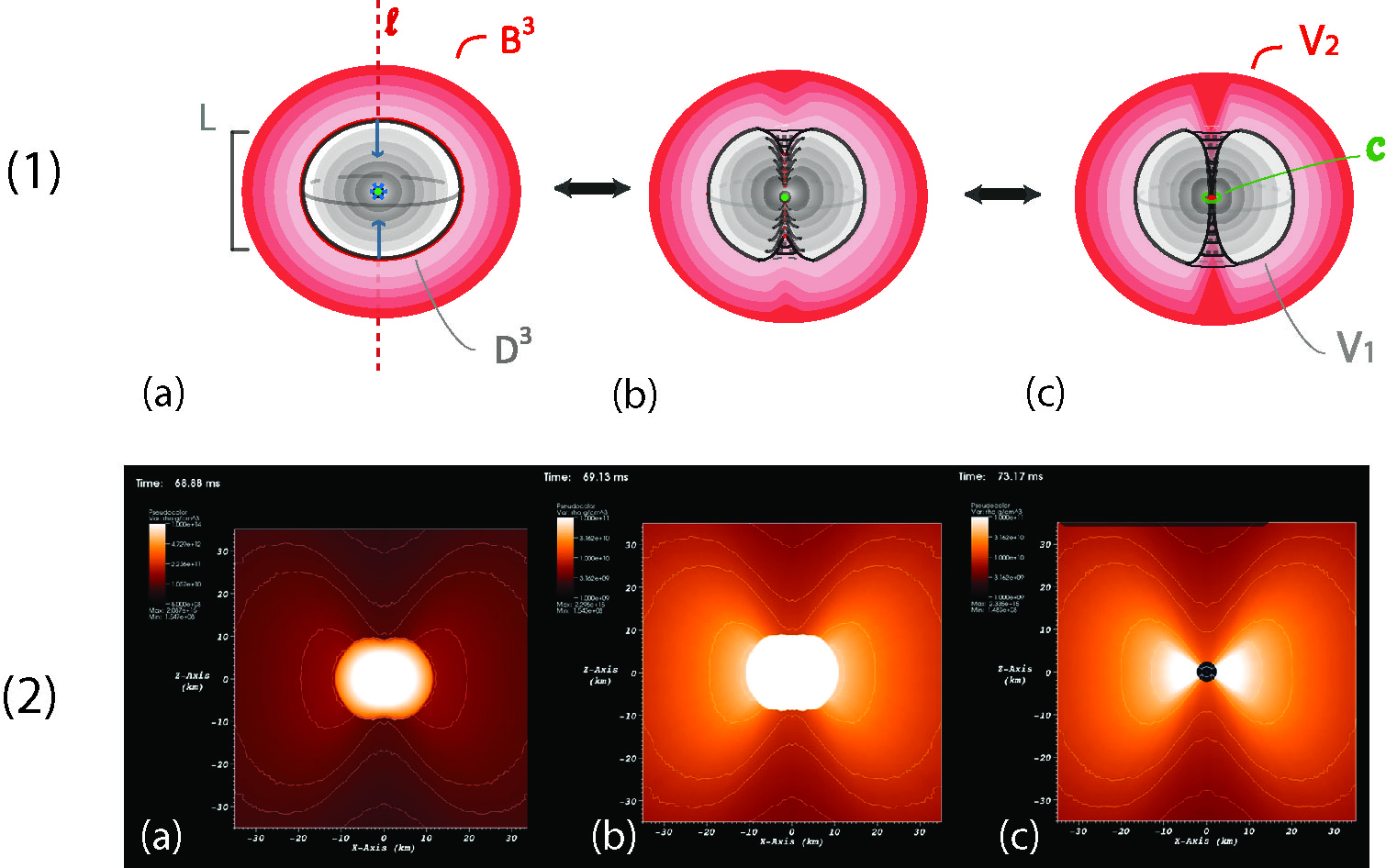}
\caption{\textbf{(1)} Embedded solid 2-dimensional 0-surgery on $M=D^3$ (in ${\mathbb R}^3$) \textbf{(2)} Black hole formation }
\label{ES2D0_B3}
\end{center}
\end{figure}

Let us start by considering the formation of {\bf black holes}. Most black holes are formed from the remnants of a large star that dies in a supernova explosion. Their gravitational field is so strong that not even light can escape. In the simulation of a black hole formation in \cite{Ott},  the density distribution at the core of a collapsing massive star is shown. Fig.~\ref{ES2D0_B3} (2) shows three  instants of this simulation,  which indicate that matter performs solid 2-dimensional 0-surgery as it collapses into a black hole. In fact,  matter collapses at the center of attraction of the initial manifold $M=D^3$ creating the singularity,  that is,  the center of the black hole (shown as a black dot in instance (c) of Fig.~\ref{ES2D0_B3} (2)),  which is surrounded by the toroidal accretion disc (shown in white in instance (c) of Fig.~\ref{ES2D0_B3} (2)). Let us be reminded here that an accretion disc is a rotating disc of matter formed by accretion.

\smallbreak

\smallbreak
\begin{figure}[!h]
\begin{center}
\includegraphics[width=12cm]{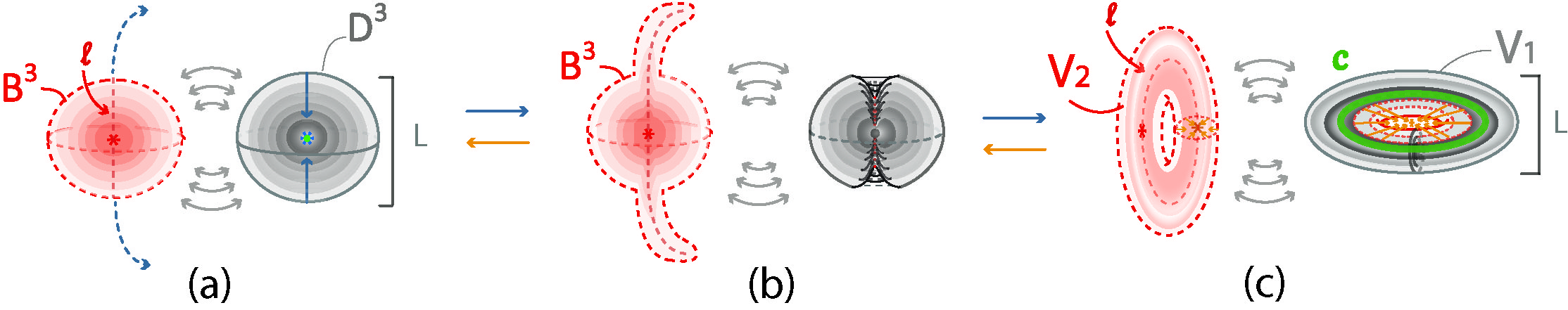}
\caption{Embedded solid 2-dimensional 0-surgery on $M=D^3$ (in $S^3$)}
\label{ES2D0_B3_Com}
\end{center}
\end{figure}

Note now that the strong gravitational forces have altered the space surrounding the initial star and that the singularity is created outside the final solid torus. This means that the process of surgery in this phenomenon has moreover altered matter outside the manifold in which it occurs. In other words,   \textit{the effect of the forces causing surgery propagates to the complement space,  thus causing a more global change in 3-space}. This fact makes black hole formation a phenomenon that topologically undergoes embedded solid 2-dimensional 0-surgery.

In Fig.~\ref{ES2D0_B3} (1),  we present a schematic model of embedded solid 2-dimensional 0-surgery on $M=D^3$. From the descriptions of $S^3$ in Section~\ref{Embedded2D},  it becomes apparent that embedded solid 2-dimensional 0-surgery on one 3-ball describes the passage from the two-ball description to the two-solid tori description of $S^3$. This can be seen in ${\mathbb R}^3$ in instances (a) to (c) of Fig.~\ref{ES2D0_B3} (1) but is more obvious by looking at instances (a) to (c) of Fig.~\ref{ES2D0_B3_Com} which show the corresponding view in $S^3$. 

We will now detail the instances of the process of embedded solid 2-dimensional 0-surgery on $M=D^3$ by referring to both the view in $S^3$ and the corresponding decompacified view in ${\mathbb R}^3$. Let $M=D^3$ be the solid ball having arc $L$ as a diameter and the complement space be the other solid ball $B^3$ containing the point at infinity; see instances  (a) of Fig.~\ref{ES2D0_B3_Com} and (a) of Fig.~\ref{ES2D0_B3}. Note that,  in both cases $B^3$ represents the hole space outside $D^3$ which means that the spherical nesting of $B^3$ in instance Fig.~\ref{ES2D0_B3}~(a) extends to infinity,  even though only a subset of $B^3$ is shown. This joining arc $L$ is seen as part of a simple closed curve $l$ passing by the point at infinity. In instances  (b) of Fig.~\ref{ES2D0_B3_Com} and (b) of Fig.~\ref{ES2D0_B3},  we see the `drilling' along $L$ as a result of the attracting forces. This is exactly the same process as in Fig.~\ref{Solid}~({b\textsubscript{1}}) if we restrict it to $D^3$. But since we have embedded the process in $S^3$ or ${\mathbb R}^3$,  the complement space $B^3$ participates in the process and,  in fact,  it is also undergoing solid 2-dimensional 0-surgery. Indeed,  the `matter' that is being drilled out from the interior of $D^3$ can be viewed as `matter' of the outer sphere $B^3$ invading $D^3$. In instances (c) of Fig.~\ref{ES2D0_B3_Com} and (c) of Fig.~\ref{ES2D0_B3},  we can see that,  as surgery transforms the solid ball $D^3$ into the solid torus $V_1$,  $B^3$ is transformed into $V_2$. That is,  the nesting of concentric spheres of $D^3$ (respectively of $B^3$) is transformed into the nesting of concentric tori in the interior of $V_1$ (respectively of $V_2$). The point at the origin (in green),  which is also the attracting center,  turns into the core curve $c$ of $V_1$ (in green) which,  by Definition~\ref{cont2D} is 2-dimensional 0-surgery on a point. As seen in instance (c) of Fig.~\ref{ES2D0_B3_Com} and (c) of Fig.~\ref{ES2D0_B3} (1),  the result of surgery is the two solid tori $V_1$ and $V_2$ forming $S^3$.

The described process can be viewed as a double surgery resulting from a single attracting center which is inside the first 3-ball $D^3$ and outside the second 3-ball $B^3$. This attracting center is illustrated (in blue) in instance (a) of Fig.~\ref{ES2D0_B3} but also in (a) of Fig.~\ref{ES2D0_B3_Com},  where it is shown that the colinear attracting forces causing the double surgery can be viewed as acting on $D^3$ (the two blue arrows) and also as acting on the complement space $B^3$ (the two dotted blue arrows),  since they are applied on the common boundary of the two 3-balls. Note that in both cases,  the attracting center coincides with the limit point of the spherical layers that $D^3$ is made of,  that is,  their common center and the center of $D^3$ (shown in green in (a) of Fig.~\ref{ES2D0_B3} and (a) of Fig.~\ref{ES2D0_B3_Com}). For more details on the descriptions of $S^3$ and their relation to surgery,  the reader is referred to \cite{SSN1}.
 
The reverse process of embedded solid 2-dimensional 0-surgery on $D^3$ is an embedded solid 2-dimensional 1-surgery on the solid torus $V_2$,  see instances of Fig.~\ref{ES2D0_B3_Com} in reverse order. This process is the embedded analog of the solid 2-dimensional 1-surgery on a solid torus $D^2 \times S^1$ defined in Definition~\ref{solto2D} and shown in Fig.~\ref{Solid}~({b\textsubscript{1}}) in reverse order. Here too,  the process can be viewed as a double surgery resulting from one attracting center which is outside the first solid torus $V_1$ and inside the second solid torus $V_2$. This attracting center is illustrated (in orange) in instance (c) of Fig.~\ref{ES2D0_B3_Com} where it is shown that the coplanar forces causing surgery are applied on the common boundary of $V_1$ and $V_2$ and can be viewed as attracting forces along a longitude when acting on $V_1$ and as attracting forces along a meridian when acting on the complement space $V_2$. 

One can now directly appreciate the correspondence of the physical phenomena (instances (a), (b), (c) of Fig.~\ref{ES2D0_B3} (2)) with our schematic model (instances (a), (b), (c) of Fig.~\ref{ES2D0_B3} (1)) . Indeed,  if one looks at the formation of black holes and examines it as an isolated event in space,  this process shows a decompactified view of the passage from a two 3-ball description of $S^3$,  that is,  the core of the star and the surrounding space,  to a two solid tori description,  namely the toroidal accretion disc surrounding the black hole (shown in white in instance (c) of Fig.~\ref{ES2D0_B3} (2)) and the surrounding space.

\begin{remark} \label{Duality2d0} \rm  It is worth pinning down the following spatial duality of embedded solid 2-dimensional 0-surgery for $M=D^3$: the attraction of two points lying on the boundary of segment $L$ by the center of $D^3$ can be equivalently viewed in the complement space as the repulsion of these points by the center of $B^3$ (that is,  the point at infinity) on the boundary of the segment $l-L$ (or the segments,  if viewed in ${\mathbb R}^3$). Hence,  the aforementioned duality tells us that the attracting forces from the attracting center that are collapsing the core of the star can be equivalently viewed as repelling forces from the point at infinity lying in the surrounding space. 
\end{remark}

\subsubsection{A topological model for the formation of tornadoes} \label{E2D0_D3S0}

\begin{figure}[!h]
\begin{center}
\includegraphics[width=11.7cm]{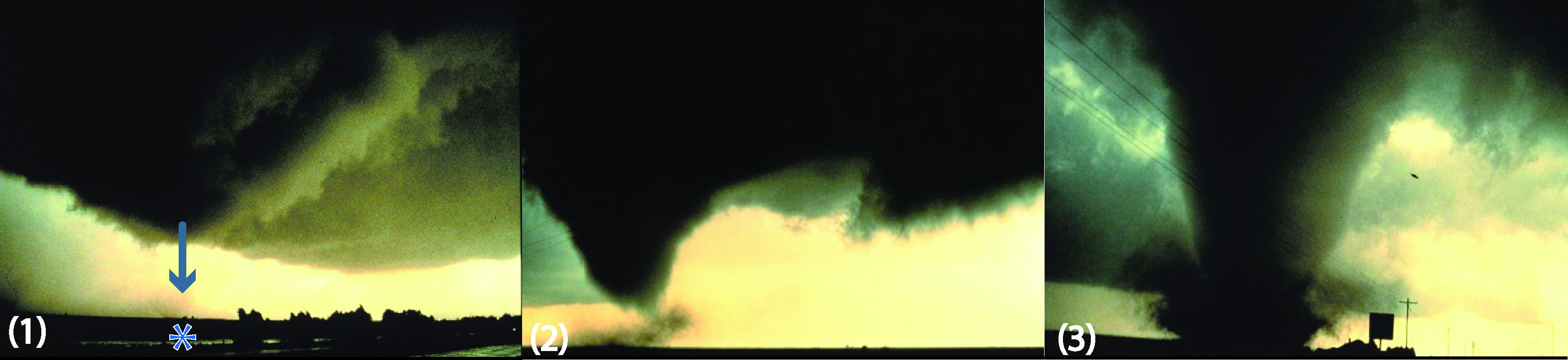}
\caption{\textbf{(1)} Attracting force between the cloud and the earth \textbf{(2)} Funnel-shaped clouds \textbf{(3)} Tornado}
\label{TornadoSequence}
\end{center}
\end{figure}

Another example of global phenomenon is the formation of {\bf tornadoes},  recall Fig.~\ref{2D_Nature} (2). As mentioned in Section~\ref{DSolid} this phenomenon can be modelled by solid 2-dimensional 0-surgery,  recall Fig.~\ref{Solid} (b\textsubscript{2}) (from right to left). However,  here,  the initial manifold is different than $D^3$. Indeed,  if we consider a 3-ball around a point of the cloud and another 3-ball around a point on the ground,  then the initial manifold is $M=D^3 \times S^0$. If certain meteorological conditions are met,  an attracting force between the cloud and the earth beneath is created. This force is shown in blue in see Fig.~\ref{TornadoSequence} (1). Then,  funnel-shaped clouds start descending toward the ground,  see Fig.~\ref{TornadoSequence} (2). Once they reach it,  they become tornadoes,  see Fig.~\ref{TornadoSequence} (3). The only difference compared to our model is that here the attracting center is on the ground,  see Fig.~\ref{TornadoSequence} (1),  and only one of the two 3-balls (the 3-ball of cloud) is deformed by the attraction. This lack of symmetry in the process can be obviously explained by the big difference in the density of the materials.

\smallbreak

During this process,  a solid cylinder $D^2\times S^1$ containing many layers of air is created. Each layer of air revolves in a helicoidal motion which is modeled using a twisting embedding as shown in Fig.~\ref{2D_Nature} (1) (for an example of a twisting embedding,  the reader is referred Section~\ref{2D_FormalE0}). Although all these layers undergo local dynamic 2-dimensional 0-surgeries which are triggered by local forces (shown in blue in Fig.~\ref{TornadoSequence} (1)),  these local forces are not enough to explain the dynamics of the phenomenon. Indeed,  the process is triggered by the difference in the conditions of the lower and upper atmosphere which create an air cycle. \textit{This air cycle lies in the complement space of the initial manifold} $M=D^3 \times S^0$ \textit{and of the solid cylinder} $D^2\times S^1$ \textit{,  but is also involved in the creation of the funnel-shaped clouds that will join the two initial 3-balls}. Therefore in this phenomenon,  surgery is the outcome of global changes and this fact makes tornado formation an example of embedded solid 2-dimensional 0-surgery on $M=D^3 \times S^0$.

\begin{figure}[!h]
\begin{center}
\includegraphics[width=11cm]{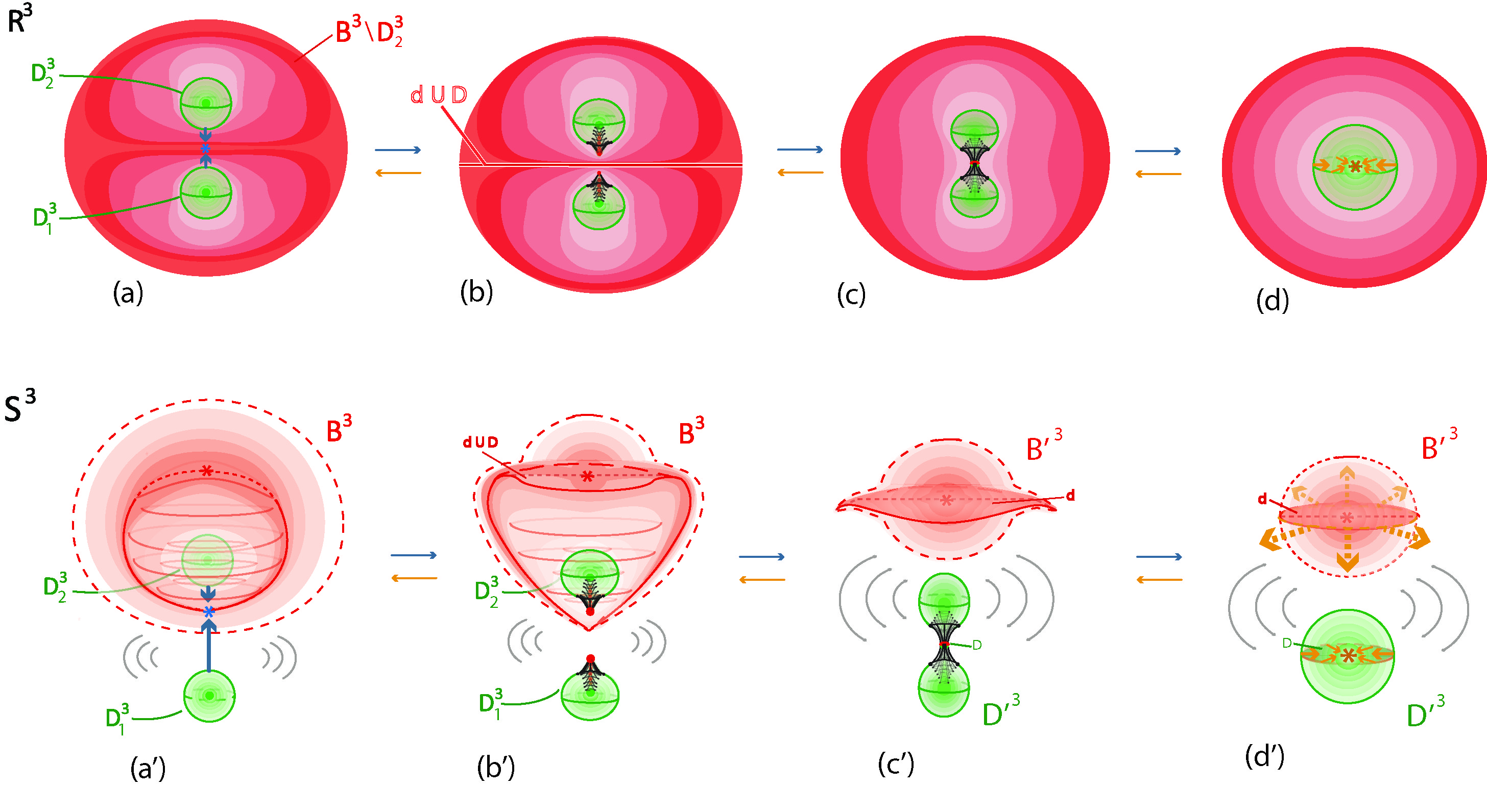}
\caption{Embedded solid 2-dimensional 0-surgery on $D^{3}_{1} \amalg D^{3}_{2}$ (from left to right) and embedded solid 2-dimensional 1-surgery on $D'^3$ (from right to left)}
\label{ES2D0_ES2D1}
\end{center}
\end{figure}

It is worth mentioning that the complement space containing the aforementioned air cycle is also undergoing solid 2-dimensional 0-surgery. The process can be seen in ${\mathbb R}^3$ in instances (a) to (d) of Fig.~\ref{ES2D0_ES2D1} while the corresponding view in $S^3$ is shown in instances (a$'$) to (d$'$) of Fig.~\ref{ES2D0_ES2D1}.


\smallbreak

More precisely,  let us name the two initial 3-balls ${D^{3}_{1}}$ and ${D^{3}_{2}}$,  hence $M=D^3 \times S^0 ={{D^{3}_{1}}} \amalg {{D^{3}_{2}}}$. Further,  let $B^3$ be the complement of ${D^{3}_{1}}$ in $S^3$. This setup is shown in (a$'$) of Fig.~\ref{ES2D0_ES2D1} where $S^3$ is viewed as the union of the two 3-balls ${D^{3}_{1}} \cup B^3$,  and here too,  $B^3$ represents everything outside ${D^{3}_{1}}$. The complement space of the initial manifold,  $S^3 \setminus M=B^3\setminus D^{3}_{2}$,  is the 3-ball $B^3$ where ${D^{3}_{2}}$ has been removed from its interior and its boundary consists in two spheres $S^2 \times S^0$,  one bounding $B^3$ or,  equivalently,  ${D^{3}_{1}}$ (the outside sphere) and one bounding ${D^{3}_{2}}$ (the inside sphere). Next,  ${D^{3}_{1}}$ and ${D^{3}_{2}}$ approach each other,  see Fig.~\ref{ES2D0_ES2D1} (b$'$). In (c$'$) of Fig.~\ref{ES2D0_ES2D1},  ${D^{3}_{1}}$ and ${D^{3}_{2}}$ merge and become the new 3-ball $D'^3$,  see Fig.~\ref{ES2D0_ES2D1} (c$'$) or  Fig.~\ref{ES2D0_ES2D1} (d$'$) for a homeomorphic representation.

At the moment of merging,  the spherical boundary of ${D^{3}_{2}}$ punctures the boundary of $B^3$; see the passage from (b$'$) to (c$'$) of Fig.~\ref{ES2D0_ES2D1}. As a result,  the complement space is transformed from $B^3\setminus D^{3}_{2}$ to the new deformed 3-ball  $B'^3$,  see Fig.~\ref{ES2D0_ES2D1} (c$'$) or Fig.~\ref{ES2D0_ES2D1} (d$'$) for a homeomorphic representation. Note that,   although the complement space undergoes a type of surgery that is different from the ones defined in Section~\ref{DSolid} and shown in Fig.~\ref{Solid},  it can still be defined analogously. In short,  we have a double solid 2-dimensional 0-surgery which turns $M={D^{3}_{1}} \amalg {D^{3}_{2}}$ into ${D'^3}$ and the complement space $S^3 \setminus ({D^{3}_{1}} \amalg {D^{3}_{2}})$ into $B'^3$. This process is initiated by the attracting center shown (in blue) in Fig.~\ref{ES2D0_ES2D1} (a$'$). The created colinear forces can be viewed as acting on ${D^{3}_{1}} \amalg {D^{3}_{2}}$ or,  equivalently,  as acting on the complement space $S^3 \setminus ({D^{3}_{1}} \amalg {D^{3}_{2}})$ (see the two blue arrows for both cases). 

Going back to the formation of tornadoes,  the above process describes what happens to the complement space and provides a topological description of the behavior of the air cycle during the formation of tornadoes. The complement space $B^3\setminus D^{3}_{2}$ in ${\mathbb R}^3$ is shown in red in Fig.~\ref{ES2D0_ES2D1} (a) and its behavior during the process can be seen in instances (b) to (d) of Fig.~\ref{ES2D0_ES2D1}. Note that in Fig.~\ref{ES2D0_ES2D1}~(a),  $B^3\setminus D^{3}_{2}$ represents the hole space outside $D^{3}_{1}$,  which means that the red layers of Fig.~\ref{ES2D0_ES2D1}~(a) extend to infinity and only a subset is shown.


\subsubsection{Embedded solid 2-dimensional 1-surgery on  $M=D^3$} \label{E2D1} 

We will now discuss the process of \textit{embedded solid 2-dimensional 1-surgery} in $S^3$. Taking $M=D'^3$ as the initial manifold,  embedded solid 2-dimensional 1-surgery is the reverse process of embedded solid 2-dimensional 0-surgery on $D^3 \times S^0$ and is illustrated in Fig.~\ref{ES2D0_ES2D1} from right to left. The process is initiated by the attracting center shown (in orange) in (d$'$) of Fig.~\ref{ES2D0_ES2D1}. The created coplanar attracting forces are applied on the circle which is the common boundary of the meridian of $D'^3$ and the meridian of $B'^3$ and they can be viewed as acting on the meridional disc $D$ of the 3-ball $D'^3$ (see orange arrows) or,  equivalently,  in the complement space,  on the meridional disc $d$ of $B'^3$ (see dotted orange arrows). As a result of these forces,  in  Fig.~\ref{ES2D0_ES2D1} (c$'$),  we see that while disc $D$ of $D'^3$ is getting squeezed,  disc $d$ of $B'^3$ is enlarged. In Fig.~\ref{ES2D0_ES2D1} (b$'$),  the central disc $d$ of $B'^3$ engulfs disc $D$ and becomes $d \cup D$,  which is a separating plane in ${\mathbb R}^3$,  see Fig.~\ref{ES2D0_ES2D1} (b). At this point the initial 3-ball $D'^3$ is split in two new 3-balls $D_1^3$ and $D_2^3$; see Fig.~\ref{ES2D0_ES2D1} (b$'$) or Fig.~\ref{ES2D0_ES2D1} (a$'$) for a homeomorphic representation. The center point of $D'^3$ (which coincides with the orange attracting center) evolves into the two centers of $D_1^3$ and $D_2^3$ (in green) which by Definition~\ref{cont2D},  is 2-dimensional 1-surgery on a point. This is exactly the same process as in Fig.~\ref{Solid}~({b\textsubscript{2}}) if we restrict it to $D'^3$,  but since we are in $S^3$,  the complement space $B'^3$ is also undergoing,  by symmetry,  solid 2-dimensional 1-surgery.


 

All natural phenomena undergoing embedded solid 2-dimensional 1-surgery take place in the ambient 3-space. The converse,  however,  is not true. For example,  the phenomena exhibiting 2-dimensional 1-surgery discussed in Section~\ref{2D} are  all embedded  in 3-space,  but they do not exhibit the intrinsic properties of embedded 2-dimensional surgery,  since they do not demonstrate the causal or consequential effects discussed in Section~\ref{Embedded2D} involving the ambient space. Yet one could,  for example,  imagine taking a solid material specimen,  stress it until necking occurs and then immerse it in some liquid until its pressure causes fracture to the specimen. In this case the complement space is the liquid and it triggers the process of surgery. Therefore,  this is an example of  embedded solid 2-dimensional 1-surgery where surgery is the outcome of global changes.

\begin{remark} \label{Duality2d1} \rm  Note that the spatial duality described in embedded solid 2-dimensional 0-surgery,  in Remark~\ref{Duality2d0},  is also present in embedded solid 2-dimensional 1-surgery. Namely,  the attracting forces from the circular boundary of the central disc $D$ to the center of $D'^3$ shown in (d$'$) of Fig.~\ref{ES2D0_ES2D1},  can be equivalently viewed in the complement space as repelling forces from the center of $B'^3$ (that is,  the point at infinity) to the boundary of the central disc $d$,  which coincides with the boundary of $D$.
\end{remark}


\begin{remark}\label{SumsUp} \rm One can sum up the processes described in this section as follows. The process of embedded solid 2-dimensional 0-surgery on $D^3$ consists in taking a solid cylinder such that the part $S^0 \times D^2$ of its boundary lies in the boundary of $D^3$,  removing it from $D^3$ and adding it to $B^3$. Similarly,  the reverse process of embedded solid 2-dimensional 1-surgery on $V_2$ consists of taking a solid cylinder such that the part $S^1 \times D^1$ of its boundary lies in the boundary of $V_2$,  removing it from $V_2$ and adding it to $V_1$. 
Following the same pattern,  embedded solid 2-dimensional 1-surgery on  $M=D^3$ consists of taking a solid cylinder in $D'^3$ such that the part $S^1 \times D^1$ of its boundary lies in the boundary of $D'^3$,  removing it from $D'^3$ and adding it to $B'^3$. Similarly,  the reverse process of embedded solid 2-dimensional 0-surgery on $S^3 \setminus ({D^{3}_{1}} \amalg {D^{3}_{2}})$ consists of taking a solid cylinder such that the two parts $S^0 \times D^2={D^{2}_{1}} \amalg {D^{2}_{2}}$ of its boundary lie in the corresponding two parts of the boundary of $S^3 \setminus ({D^{3}_{1}} \amalg {D^{3}_{2}})$,  removing it from $S^3 \setminus ({D^{3}_{1}} \amalg {D^{3}_{2}})$ and adding it to ${D^{3}_{1}} \amalg {D^{3}_{2}}$. 
Note that,  for clarity,  in the above descriptions the attracting centers causing surgery are always inside the initial manifold. Of course a similar description starting with the complement space as an initial manifold and the attracting center outside of it would also have been correct.
\end{remark}

\section{Conclusions}
In this paper,  inspired by natural phenomena exhibiting topological surgery,  we introduced the notions of dynamic,  solid and embedded surgery in dimensions 1,  2 and 3. There are many more natural phenomena exhibiting surgery and we believe that the understanding of their underpinning topology will lead to the better understanding of the phenomena themselves,  as well as to new mathematical notions,  which will in turn lead to new physical implications.




\begin{acknowledgement}  
We wish to express our gratitude to Louis H.Kauffman and Cameron McA.Gordon for many fruitful conversations on topological surgery. We would also like to thank the Referee for his/her positive comments and for helping us clarify some key notions.  \\
We further wish to acknowledge that this research  has been co-financed by the European Union (European Social Fund - ESF) and the Greek national funds through the Operational Program "Education and Lifelong Learning" of the National Strategic Reference Framework (NSRF) - Research Funding Program: THALES: Reinforcement of the interdisciplinary and/or inter-institutional research and innovation.
\end{acknowledgement}



\end{document}